\documentclass[oneside,a4paper,12pt,reqno]{amsart}
\usepackage{upref,amsfonts,amsxtra,amssymb}

\newcommand{\Z}{\mathbb{Z}}

\newtheorem{theorem}{Theorem}[section]
\newtheorem{lemma}[theorem]{Lemma}
\newtheorem{prop}[theorem]{Proposition}

\newtheorem{q}[theorem]{Question}
\title{Periodic points of algebraic actions of discrete groups}
\author{Siddhartha Bhattacharya}
\address{School of Mathhematics, Tata Institute of Fundamental
  Research, Mumbai 400005, India}
\email{siddhart@math.tifr.res.in}
\subjclass[2010]{37B05, 37B20}
\keywords{Nilpotent groups, automorphisms, periodic points}
\date{}
\begin{document}
\maketitle
\begin{abstract}
Let $\Gamma$ be a countable group. A $\Gamma$-action on a 
compact abelian group $X$ by continuous automorphisms of
$X$ is called Noetherian if the dual of $X$ is Noetherian 
as a  ${\mathbb Z}(\Gamma)$-module. 
We prove that any Noetherian action of a finitely generated
virtually nilpotent group has a dense set of periodic points.
\end{abstract}
\section{Introduction}
Let $\Gamma$ be a countable group. An {\it algebraic $\Gamma$-action\/} 
is an action $\alpha$ 
of $\Gamma$ on a compact metrizable abelian group $X$
by continuous automorphisms of $X$. By duality theory, any such action
induces a $\Gamma$-action on ${\widehat X}$, the dual of $X$, by
automorphisms of ${\widehat X}$. Hence ${\widehat X}$ can be viewed as
a $\Z[\Gamma]$-module, where $\Z[\Gamma]$ is the integral group ring
of $\Gamma$. The action $\alpha$ is called {\it  Noetherian\/} if
${\widehat X}$ is Noetherian as a $\Z[\Gamma]$-module. Equivalently, 
$\alpha$ is Noetherian if any decreasing sequence
$$ X = X_0 \supset X_1\supset X_2\supset\cdots $$
of closed $\Gamma$-invariant subgroups stabilizes.

For $\Gamma = \Z^d$, the study of such actions was initiated in
\cite{KS}, and since then these systems have been extensively studied
(see \cite{Sc} for a comprehensive account). For general $\Gamma$,
several important results about entropy and homoclinic points of 
Noetherian $\Gamma$-actions have been proved in recent years
(\cite{Bo}, \cite{CL}, \cite{De}, \cite{Li}). 
However, unlike the commutative case, many
basic dynamical properties are only partly understood in the general
situation.
 
In this paper, we study periodic points of Noetherian
$\Gamma$-actions. Recall that a point $x\in X$ is periodic if the
$\Gamma$-orbit of $x$ is finite. It is known that when $\Gamma = \Z^d$, for any
Noetherian $\Gamma$-action $(X,\alpha)$ the space $X$ contains a dense
set of periodic points (\cite{Sc}, Theorem 5.7). However, the proof
uses tools from commutative algebra and does not generalize to actions
of non-abelian groups.
The question whether Noetherian actions of all residually finite groups
admit a dense set of periodic points was raised in
\cite[Problem 8.5]{LS}. 

In this paper we show that for algebraic actions of countable groups,
the density of periodic points is related to vanishing of certain
cohomology groups. As a consequence, we obtain the following : 
\begin{theorem}\label{nil}
Let $\Gamma$ be a finitely generated virtually nilpotent group and let 
$\alpha$ be a Noetherian action of $\Gamma$ on a compact metrizable
abelian group $X$. Then $X$ contains a dense set of $\alpha$-periodic points. 
\end{theorem}
As another application we show that there exist Noetherian actions of
finitely generated residually finite groups that do not have
a dense set of periodic points, thus giving a negative answer to the
above question in the most general case.

\section{Virtual first cohomology groups}
Let $\alpha$ be an action of a countable group $\Gamma$ on a compact
metrizable group $X$ by continuous automorphisms of $X$. A 
{\it virtual 1-cocycle \/} of $\alpha$ is a map $c$ from a
finite index subgroup $\Lambda\subset \Gamma$ to $X$ that
satisfies the equation 
$$ c(\gamma \gamma^{'}) = c(\gamma) + \alpha(\gamma )(c(\gamma^{'}))$$
for all $\gamma, \gamma^{'}\in \Lambda$. Two 
virtual 1-cocycles $c_1 : \Lambda_1\rightarrow X$ and $c_2 :
\Lambda_2\rightarrow X$ are {\it equivalent\/} if there exists a finite index
subgroup $\Lambda\subset  \Lambda_1\cap \Lambda_2$ such that  
$c_1(\gamma) = c_2(\gamma)$ for all $\gamma\in \Lambda$. 
For any $x\in X$, the map $c_x : \Gamma\rightarrow X$ defined by 
$c_x(\gamma) = \alpha(\gamma)(x) - x$ is a 1-cocycle. 
A virtual 1-cocycle $c : \Lambda\rightarrow X$ 
is said to be a {\it virtual coboundary \/} if it is equivalent to
$c_x$ for some $x\in X$. We will call 
two virtual 1-cocycles $c_1 : \Lambda_1\rightarrow X$ and $c_2 :
\Lambda_2\rightarrow X$  {\it cohomologous\/} if there exists a finite 
index subgroup $\Lambda\subset \Lambda_1\cap \Lambda_2$, and a virtual 
coboundary $c : \Lambda\rightarrow X$ such that $c_1(\gamma) -
c_2(\gamma) = c(\gamma)$ for all $\gamma\in\Lambda$.
 It is easy to see that the equivalence classes
of virtual 1-cocyles is a group with respect to the pointwise
addition, and the equivalence classes of virtual coboundaries form a subgroup. 
We will denote the quotient group 
by $H_{v}^{1}(\alpha)$ and call it the {\it virtual 
first cohomology group \/} of $\alpha$.

Let $(Y, \beta)$ be an algebraic $\Gamma$-action,
and let $X\subset Y$ be a closed $\beta$-invariant
subgroup such that $Y/X$ is finite. Let $\alpha$ denote the
restriction of $\beta$ to $X$ and let $\Lambda\subset \Gamma$ be 
a finite index subgroup that acts trivially on $Y/X$. 
For any $y\in Y$ we define a map $c_y : \Lambda
\rightarrow X$ by $c_y(\gamma) = \beta(\gamma)(y) - y$. It is easy to
see that $c_y$ is a virtual 1-cocycle of $\alpha$. 

\begin{prop}\label{extco}
Let $(Y,\beta)$, $(X,\alpha)$ and $y\in Y$ be as above. Then the
virtual 1-cocycle $c_y$ is a virtual coboundary if and only if the
coset $y + X$ contains a $\beta$-periodic point. 
\end{prop}

{\it Proof.\/} Suppose $y_1\in y + X$ is a $\beta$-periodic point, and
$\Lambda_1\subset \Lambda$ is a finite index subgroup such that 
$\beta(\gamma)(y_1) = y_1$ for all $\gamma\in\Lambda_1$. Then $x = y -
y_1\in X$ and 
$$c_y(\gamma) = \beta(\gamma)(y) - y = \beta(\gamma)(x) - x$$
for all $\gamma\in\Lambda_1$. Hence $c_y$ is a virtual coboundary of $\alpha$. 
Conversely, if $c_y$ is a virtual coboundary of $\alpha$,
 then there exists 
a finite index subgroup $\Lambda_1$ and $x\in X$ such that 
$c_y(\gamma) = \beta(\gamma)(x) - x$ for all $\gamma\in\Lambda_1$. If $p
= y - x$ then $p\in y + X$ and 
it is fixed by $\beta(\Lambda_1)$. Hence it is a periodic point. $\Box$

\medskip
\noindent
A compact abelian group $X$ is {\it zero-dimensional \/} if the
connected component containing $0_X$ is trivial. From the duality
theory it follows that $X$ is zero-dimensional if and only if 
every element of ${\widehat X}$ has finite order. 
\begin{prop}\label{torsion}
Let $\Gamma$ be a countable group, and let $\alpha$ be a Noetherian
action of $\Gamma$ on a zero-dimensional compact abelian group $X$.
Then there exists $k\ge 1$ such that $kx = 0$ for all $x\in X$.
\end{prop}
{\it Proof.\/}
 Let $\{ \chi_1,\ldots, \chi_n\}\subset {\widehat X}$ be a finite set 
that generates ${\widehat X}$ as a ${\mathbb Z}(\Gamma)$-module. As
every element of ${\widehat X}$ has finite order, we can find
 $k\ge 1$
such that $k\chi_i = 0$ for $i = 1,\ldots , n$. Since for any $m$ the set 
$N_m = \{ \chi : m\chi = 0\}$ is a $\Z(\Gamma)$-submodule of ${\widehat X}$, we
deduce that $N_k = {\widehat X}$, i.e., $k\chi = 0$ for all $\chi\in
{\widehat X}$. If $x\in X$ then for all 
$\chi \in {\widehat X}$, $\chi(kx) = k\chi(x) = 0$. Hence $kx = 0$. 
$\Box$ 

\medskip
\noindent
Our next proposition shows that all virtual 1-cocycles of Noetherian
$\Gamma$-actions on zero-dimensional groups arise in the manner
described in Proposition \ref{extco}.
\begin{prop}\label{coext}
Let $\Gamma$ be a countable group and let $(X,\alpha)$ be a Noetherian
action of $\Gamma$ on a zero-dimensional group $X$. Then for any
virtual 1-cocycle $ c : \Lambda\rightarrow X$ there exists $k\ge 1$ and an 
algebraic $\Lambda$-action $\beta$ on $Y = ({\mathbb Z}/k{\mathbb
    Z})\times X$ such that
\begin{enumerate}
\item{For all $\gamma\in\Lambda$, $\beta(\gamma)$ fixes the first co-ordinate.}
\item{$\beta(\gamma)(0,x) = (0, \alpha(\gamma)(x))$ for all $\gamma\in\Lambda$
    and $x\in X$.}
\item{$c_{(1,0)} = c.$}

\end{enumerate}
\end{prop}

{\it Proof.\/} Since $\alpha$ is Noetherian and $X$ is zero-dimensional,
by the previous proposition there exists $k\ge 1$ such that $kx = 0$ for all
$x\in X$. For any $z\in Y$ and 
a continuous automorphism $\tau$ of $X$ there exists
  a unique homomorphism $\theta : Y \rightarrow Y$ such
  that $\theta(1,0) = z$ and $\theta(0,x) = (0,\tau(x))$ for all $x\in
  X$. It is easy to see that $\theta$ is continuous.
 For $\gamma\in\Lambda$ let $\beta(\gamma)$ denote the unique
 continuous endomorphism of
 $Y$ satisfying the following two conditions :
\begin{enumerate}
\item{$\beta(\gamma)(0,x) = (0, \alpha(\gamma)(x))$ for all $x\in X$.}
\item{$\beta(\gamma)(1,0) = (1, c(\gamma))$.}
\end{enumerate}
Similarly we define $\beta^{'}(\gamma)$ to be the unique continuous
endomorphism of $Y$ with the property that 
$\beta^{'}(\gamma)(0,x) = (0, \alpha(\gamma)^{-1}(x))$ and 
$\beta^{'}(\gamma)(1,0) = (1, -\alpha(\gamma)^{-1}c(\gamma))$.
Then 
$$\beta(\gamma)\beta^{'}(\gamma)(1,0) = \beta(\gamma)(1,0)
+ \beta(\gamma)(0, -\alpha(\gamma)^{-1}c(\gamma)) = (1,0).$$
As $\beta(\gamma)\beta^{'}(\gamma)$ fixes both $\{ 0\} \times X$ and
the point $(1,0)$, we conclude that $\beta(\gamma)\beta^{'}(\gamma) =
I$. A similar computation shows that $\beta^{'}(\gamma)\beta(\gamma) =
I$. Hence $\beta(\gamma)$ is a continuous automorphism of $Y$ 
for each $\gamma\in\Lambda$. As 
$c(\gamma_1\gamma_2) = c(\gamma_1) + \alpha(\gamma_1)c(\gamma_2)$, it
follows that 
$$\beta(\gamma_1\gamma_2) (1,0) = (1, c(\gamma_1\gamma_2) ) = 
\beta(\gamma_1)\beta(\gamma_2)(1,0).$$
Since $\beta(\gamma_1\gamma_2)(0,x) =
\beta(\gamma_1)\beta(\gamma_2)(0,x)$ for all $x\in X$, this imlies
that $\beta(\gamma_1\gamma_2) = \beta(\gamma_1)\beta(\gamma_2)$.
Hence
$\beta$ defines an algebraic $\Lambda$-action on $Y$. It is easy to see
that $\beta$ satisfies all three conditions stated in the
proposition. $\Box$

\begin{lemma}\label{periodic}
Let $\Gamma$ be a countable group and let $(X,\alpha)$ be a Noetherian
$\Gamma$-action such that $X$ admits non-zero periodic
points. Then $X$ also contains non-trivial $\alpha$-invariant finite
subgroups.
\end{lemma}

{\it Proof.\/}
Let $x\in X$ be a non-zero $\alpha$-periodic point and let
$\Lambda\subset \Gamma$ be the stabilizer of $x$. Since the orbit of
$x$ is finite, $\Lambda $ is a finite index subgroup of $\Gamma$. We
define a subgroup $\Lambda_0\subset \Gamma$ by 
$$\Lambda_0 =  \bigcap_{\gamma\in \Gamma}\gamma\Lambda \gamma^{-1}.$$
It is easy to see that $\Lambda_0\subset \Gamma$ is normal and has
finite index. Let $Y$ denote the set of points of $X$ that are fixed
by all elements of $\Lambda_0$. Then $Y$ is a non-trivial $\alpha$-invariant 
closed subgroup of
$X$. Since $\Lambda_0\subset \Gamma$ has finite index,
the action $\alpha|_{\Lambda_0}$ is also Noetherian. As $\Lambda_0$ acts
trivially on $Y$, it follows that any collection of closed subgroups
of $Y$ has a minimal element. Let $Y_0$ be a minimal element of the
collection of all non-trivial subgroups of $Y$. Since closed subgroups
of $Y_0$ are in one to one correspondence with subgroups of ${\widehat
  Y_0}$, we deduce that ${\widehat Y_0}$ does not admit non-trivial
subgroups. Therefore ${\widehat Y_0}$, and hence $Y_0$, are isomorphic with 
${\mathbb Z}/p{\mathbb Z}$ for some prime $p$. In particular, $Y_0$ is
finite. We define $H\subset X$ by 
$$H = \sum_{\gamma\in\Gamma}\alpha(\gamma)(Y_0).$$
As $\alpha(\gamma)(Y_0) = Y_0$ for all $\gamma\in \Lambda_0$,
and $\Lambda_0\subset\Gamma$ has finite index, the above sum is finite
and $H$ is a well defined
closed subgroup of $X$. It is easy to see that $H$ is finite and
$\alpha$-invariant. $\Box$
  
\medskip
\noindent
The above lemma is not true if $\alpha$ is not assumed to be
Noetherian. For example, suppose $\Gamma$ is an arbitrary countable
group  and $\alpha$ is the trivial action of $\Gamma$ on 
$X = {\widehat {\mathbb Q}}$. Then every point of $X$ is
$\alpha$-periodic, but since the group ${\mathbb Q}$ does not have
non-trivial finite quotients, $X$ does not admit non-trivial finite
subgroups. 

\begin{theorem}\label{ns}
Let $\Gamma$ be a countable group satisfying :
\begin{enumerate}
\item{ Every Noetherian $\Gamma$-action $(X,\alpha)$ with $X\ne \{
0\}$, has a non-zero $\alpha$-periodic point.}
\item{The group $H_{v}^{1}(\alpha)$ is trivial for every Noetherian 
$\Gamma$-action $(X,\alpha)$ on a zero-dimensional group $X$.}
\end{enumerate}
Then every Noetherian action of $\Gamma$ on a compact abelian group
(not necessarily zero-dimensional) admits a dense set of periodic
points. Conversely, if every Noetherian $\Gamma$-action
admits a dense set of periodic points, then both the conditions are satisfied. 
\end{theorem}
{\it Proof.\/} Suppose every Noetherian action of 
$\Gamma$ admits a dense set of periodic orbits. Then the first
condition is automatically satisfied. Let $\alpha$ be a 
Noetherian $\Gamma$-action  on a zero-dimensional group $X$, and let 
$c : \Lambda\rightarrow X$ be a virtual 1-cocycle of $\alpha$. By Proposition
\ref{coext}, there exists $k\ge 1$ and an algebraic $\Lambda$-action $\beta$ on 
$Y = ({\mathbb Z}/k{\mathbb Z})\times X$ such that $c = c_{(1,0)}$ and 
$\beta(\gamma)|_{X} = \alpha(\gamma)$ for all $\gamma\in\Lambda$.  
Since the coset containing $(1,0)$ is an open subset of $Y$, it
contains a $\beta$-periodic orbit. From Proposition \ref{extco}  we deduce that
$c$ is a virtual coboundary. Hence $H^{1}_{v}(\alpha) = \{0\}$. 

Now suppose $\Gamma$ satisfies both the conditions stated above. 
Let $(Y,\beta)$ be a Noetherian action of $\Gamma$. 
Let $X\subset Y$ denote the closure of the set
of all $\beta$-periodic points. 
Since the sum of two periodic points is again a
periodic point
it follows that  $X$ is a closed $\beta$-invariant
subgroup. Suppose $X\ne Y$. Let ${\bar \beta}$ denote the induced
action of $\Gamma$ on the quotient $Y/X$. 
Since $\Gamma$ satisfies the first condition, by Lemma \ref{periodic} there 
exists a finite ${\bar \beta}$-invariant non-zero subgroup  
$F\subset Y/X$. We define $Y^{'} = \pi^{-1}(F)$, where $\pi$ is the
projection map from $Y$ to $Y/X$.
We choose a finite index subgroup 
$\Lambda\subset \Gamma$ that acts trivially on $F$ under ${\bar \beta}$. 
Let $y\notin X$ be a point in $Y^{'}$, and let $c_y : \Lambda\rightarrow
X$ denote the virtual 1-cocycle of $(X,\beta)$ defined by 
$c_y(\gamma) = \beta(\gamma)(y) - y$.

Let $X_0$ denote the connected component of $X$. It is easy to see
that $X_0$ is a $\beta$-invariant closed subgroup of $X$. Let
$\beta_1$ denote the induced action of $\Lambda$ on $X/X_0$.
As $X/X_0$ is zero-dimensional and $\beta_1$ is Noetherian, applying
Proposition \ref{torsion}
 we deduce that there exists $l\ge 2$ such that $lp \in X_0$
for all $p\in X$. Let $k\ge2$ be a positive integer such that $kx = 0$
for all $x\in F$. Then $lky\in X_0$. Since $X_0$ is connected, the map 
$x\mapsto lkx$ is a surjective endomorphism of $X_0$. 
We find $q\in X_0$ such that $lky =
lkq$, and define a virtual 1-cocycle $c_1 : \Lambda \rightarrow X$ by
$c_1 = c_y - c_q$. Then for any $\gamma\in \Lambda$,
$$ lk c_1(\gamma) = lk(\beta(\gamma)(y - q) - (y-q)) =0.$$
Hence the image of $c_1$ is contained in $ H = \{ x\in X : lkx = 0\}$.
Since $H$ is closed, $\beta$-invariant, and zero-dimensional; from the
second condition we conclude that $c_1$ is a
virtual coboundary. As $c_1$ is cohomologous to $c_y$, this implies
that $c_y$ is also a virtual coboundary. Applying Proposition \ref{extco} we
deduce that the coset $y + X$ contains a periodic point
$y_1$. This contradicts the fact that $y\notin
X$, and shows that $Y = X$, i.e., the set of $\beta$-periodic points
is dense in $Y$.
$\Box$

\medskip
We now construct an example of a Noetherian action of a finitely
generated residually finite group that does not have a dense set of
periodic points. 
Let $H$ denote the group of all functions from ${\mathbb Z}$ to
${\mathbb Z}$ with finite support, equipped with pointwise
addition. We define an action of ${\mathbb Z}$ on $H$ by
$ n\cdot f(i)  = f(i + n).$ Let $\Gamma$ denote the semi-direct
product of ${\mathbb Z}$ and $H$ defined by
$$(f_1, n_1) \cdot (f_2, n_2) = ( f_1 + n_1\cdot f_2, n_1 + n_2).$$
It is easy to see that $\Gamma$ is torsion free.
For $k\in {\mathbb Z}$, we define $f^k\in H$ by $f^k(i) = 0$ if 
$k\ne i$ and $f^k(i) = 1$ if $k = i$. Then for any $k\in {\mathbb Z}$,
$(0,1)(f^{k+1},0) =  (f^{k},0)(0,1)$, i.e.,
$(0,1)(f^{k+1},0)((0,1)^{-1} = (f^{k},0)$.
As $\{ (0, f^k) : k\in {\mathbb Z}\}$ generates $H$ as a ${\mathbb
  Z}$-module, this shows that $\Gamma$ is generated by $(0,1)$ and
$(f^{0},0)$.
In particular, $\Gamma$ is finitely generated.
For $k\ge 2$ let $\Gamma_k\subset \Gamma$ denote the set of all
elements $(f, n)$ such that $n = 0 $ (mod $k$) and
$$\sum_{i= -\infty}^{\infty}f(ki + j) = 0 ({\rm mod }k)\ \forall j =
0,1, \ldots, k-1.$$
It is easy to see that for each $k$, $\Gamma_k$ is a 
normal subgroup of $\Gamma$ and $\Gamma/\Gamma_k$ is finite. 
Moreover for each non-zero $(f,n)\in \Gamma$ there exists $k$ such
that $(f,n)\notin \Gamma_k$. Hence $\Gamma$ is residually finite. 

Let $X$ denote the compact abelian group 
$({\mathbb Z}/2{\mathbb Z})^{{\mathbb Z}}$, equipped with the 
product topology and pointwise
addition, and let $S : X\rightarrow X$ denote the shift map defined
by $S(x)(i) = x(i+ 1)$. We define an algebraic $\Gamma$-action
$\alpha$ on $X$ by $\alpha(f,n)(x) = S^{n}(x)$. It is easy to see that 
the shift action of $\Z$ on $X$ is Noetherian. Since $\alpha|_{H}$ is
trivial and the action of $\Gamma/H \cong \Z$  induced by $\alpha$ is
the shift action on $X$, we deduce that $\alpha$ is also Noetherian.

 Let $\pi : H\rightarrow X$ denote the homomorphism defined by
$\pi(f)(i) = f(i)\  (mod\ 2)$. We define $c : \Gamma\rightarrow X$ by 
$c(f, n) = \pi(f)$. Then
$c((f_1, n_1)(f_2, n_2)) = \pi(f_1 + n_1\cdot f_2) = c(f_1,n_1) +
S^{n_1}(\pi(f_2))$.
Since 
$$S^{n_1}(\pi(f_2)) = S^{n_1}(c(f_2,n_2)) = \alpha(f_1, n_1)(c(f_2,
n_2)),$$
this shows that $c$ is a 1-cocycle of $\alpha$. 
We pick any $x\in X$ and define $c^{'} : \Gamma \rightarrow X$ by
$c^{'}(\gamma) = \alpha(\gamma)(x) - x$. Let $\Lambda$ be an arbitrary
finite index subgroup of $\Gamma$. Since $\alpha(\gamma) = I$ for all
$\gamma\in H$, it follows that $c^{'}(H\cap\Lambda) = \{ 0\}$. On the
other hand, the restriction of $c$ to $H$ is a homomorphism with
infinite image. As $H\cap \Lambda$ is a
finite index subgroup of $H$, we deduce that  $c(H\cap\Lambda)$ is also
infinite. Hence $c|_{\Lambda} \ne c^{'}|_{\Lambda}$. As $\Lambda$ and
$x$ are arbitrary, we conclude that $c$ is not a virtual
coboundary. We now apply Proposition \ref{coext}. Let
 $\beta$ denote the action corresponding to the cocycle
$c$. Since $\alpha$ is Noetherian, so is $\beta$.
Applying  Proposition
\ref{extco} we conclude that $\beta$
does not have a dense set of periodic points.

\section{Density of periodic orbits}

In this section we concentrate on the case when $\Gamma$ is 
polycyclic-by-finite. 
A countable group $\Gamma$ is {\it polycyclic\/} if there exists a
decreasing sequence of subgroups
$$ \Gamma = \Gamma_n\supset \Gamma_{n-1}\supset \cdots \supset
\Gamma_0 = \{ 0\}$$
such that for each $i$, 
$\Gamma_i$ is a normal subgroup of $\Gamma_{i+1}$ and
$\Gamma_{i+1}/\Gamma_{i}$ is cyclic. Any such series is  called a {\it
  polycyclic series\/} of $\Gamma$. A group $\Gamma$ is called
{\it polycyclic-by-finite\/}
if it contains a finite index subgroup that 
is polycyclic. If $\Gamma$ is polycyclic-by-finite then 
every subgroup of $\Gamma$ is finitely generated, and
${\mathbb Z}(\Gamma)$ is a Noetherian ring.

Let $\Gamma$ be a polycyclic-by-finite group, $\Gamma_0\subset
\Gamma$ be a polycyclic subgroup of finite index and
$ \Gamma_0 = \Gamma_n\supset \Gamma_{n-1}\supset \cdots \supset
\Gamma_0 = \{ 0\}$ be a polycyclic series of $\Gamma_0$. Then
 the number of $i$'s such that $\Gamma_{i+1}/\Gamma_i$ is
infinite cyclic is independent of $\Gamma_0$ and the polycyclic
series. This number is known as the {\it Hirsch number \/} of
$\Gamma$. The following proposition summarizes some basic properties
of this invariant (see \cite{W}).
  
\begin{prop}\label{Hirsch}
Let $\Gamma$ be a polycyclic-by-finite group. 
\begin{enumerate}
\item{$h(\Gamma) = 0$ if and only if $\Gamma$ is finite.}
\item{If $\Gamma_1$ is a subgroup of $\Gamma$ then $h(\Gamma_1)\le
    h(\Gamma)$.}
\item{If $\Gamma_1\subset \Gamma$ is  normal then 
$h(\Gamma) = 
    h(\Gamma_1) + h(\Gamma/\Gamma_1)$.}
\end{enumerate}
\end{prop}
We note that the previous proposition  applies to
 finitely generated virtually nilpotent groups since they are
 polycyclic-by-finite. In the proof of Theorem 1.1 we will also use the 
following result about polycyclic-by-finite groups (see \cite{Ro}): 
\begin{prop}\label{module}
Let $\Gamma$ be a polycyclic-by-finite group and let $M$ be a simple 
$\Z[\Gamma]$-module. Then $M$ is finite.
\end{prop}

Our next lemma is a direct consequence of this result.

\begin{lemma}\label{simple}
Let $\Gamma$ be a polycyclic-by-finite group and let $(X,\alpha)$ be a
Noetherian action of $\Gamma$ on a non-trivial 
compact abelian group $X$. Then $X$
admits a non-zero periodic point.
\end{lemma}

{\it Proof.\/} Let ${\mathcal A}$ be the collection of all proper ${\mathbb
  Z}(\Gamma)$-submodules of ${\widehat X}$. Since $\alpha$ is
Noetherian, ${\mathcal A}$ contains a maximal element $M$. It is easy to
see that $N = {\widehat X}/M$ is a simple ${\mathbb
  Z}(\Gamma)$-module. By the previous proposition $N$ is finite. Let 
$i : {\widehat N}\rightarrow X$ denote the dual of the projection map
$\pi : {\widehat X}\rightarrow N$. Since $\pi$ is surjective, $i$ is
injective. Hence $i({\widehat N})$ is a non-trivial $\alpha$-invariant 
finite subgroup of $X$. Clearly every non-zero point of $i({\widehat
  N})$ is periodic. $\Box$

\medskip
\noindent
\begin{lemma}\label{lem}
Let $\Gamma$ be a polycyclic-by-finite
 group and let $\alpha$ be a Noetherian action of 
$\Gamma$ on a zero-dimensional group $X$ such that 
$H^{1}_{v}(X, \alpha) \ne \{ 0\}$. If ${\rm Ker}(\alpha) = \{ \gamma :
\alpha(\gamma) = I\}$ is infinite then there exists a
polycyclic-by-finite group $\Gamma^{'}$ with $h(\Gamma^{'}) < h(\Gamma)$, and a
Noetherian $\Gamma^{'}$-action $\alpha^{'}$ on $X$ such that
$H^{1}_{v}(X, \alpha^{'}) \ne \{ 0\}$.
\end{lemma}

{\it Proof.\/} Let $\Lambda\subset\Gamma$ be a finite index subgroup,
and let $c :\Lambda\rightarrow X$ be a virtual 1-cocycle of $\alpha$
that is not a virtual coboundary. Replacing $\Gamma$ by $\Lambda$ if
necessary, we may assume that $\Lambda = \Gamma$.
 Since ${\rm Ker }(\alpha)$ acts
trivially on $X$, it follows that the restriction of $c$ to ${\rm Ker
}(\alpha)$ is a homomorphism. Since $\Gamma$ is polycyclic-by-finite,
the subgroup 
$ {\rm Ker }(\alpha)$ is finitely generated. This implies that 
$c({\rm Ker}(\alpha))$ is also finitely generated.
By Proposition \ref{torsion} there
exists $k\ge 2$ such that $kx = 0$ for all $x\in X$. Hence $c({\rm
  Ker}(\alpha))$ is a finite subgroup of $X$. We define 
$M = \{ \gamma\in {\rm Ker}(\alpha) : c(\gamma) = 0\}$. 
 Since $ {\rm Ker }(\alpha)$ is finitely generated and $M\subset {\rm
   Ker}(\alpha)$ has finite index, $M$ contains a finite
index characteristic subgroup of $ {\rm Ker }(\alpha)$, 
i.e., there exists a finite index
subgroup $N\subset M$ such that $\theta(N) = N$ for
 all automorphism $\theta$ of ${\rm Ker }(\alpha)$. 
Since for any $g\in \Gamma$ the map 
$\gamma\mapsto g\gamma g^{-1}$ is an automorphism of ${\rm
  Ker}(\alpha)$, we deduce that $N$ is a normal subgroup of $\Gamma$. 
We define $\Gamma^{'} = \Gamma/N$. As $N$ is a finite index subgroup
of ${\rm Ker}(\alpha)$ and ${\rm Ker}(\alpha)$ is infinite, $N$ is
also infinite. In particular $h(N) > 0$. Hence 
$$h(\Gamma/N) = h(\Gamma) - h(N) < h(\Gamma).$$
As $\alpha(\gamma) = I$ for all $\gamma\in N$, $\alpha$ induces a
Noetherian $\Gamma^{'}$-action $\alpha^{'}$ on $X$. We define 
$c^{'}: \Gamma^{'}\rightarrow X$ by $c^{'}(\gamma N) = c(\gamma)$.
We note that for any $\gamma\in\Gamma$ and $n\in N$, 
$c(\gamma n) = c(\gamma) + \alpha(\gamma)(c(n)) = c(\gamma)$. Hence
$c^{'}$ is a well defined virtual 1-cocycle of $\alpha^{'}$. Suppose 
there exist a finite index subgroup $\Lambda^{'}\subset \Gamma^{'}$ and
$x\in X$ such that $c^{'}(a) = \alpha^{'}(a)(x) - x$ for all $a\in
\Lambda^{'}$. We define $\Lambda_1 = \pi^{-1}(\Lambda^{'})$, where 
$\pi : \Gamma\rightarrow \Gamma^{'}$ is the projection map. Then
for any $\gamma\in\Lambda_1$,
$$c(\gamma) = c^{'}(\pi(\gamma)) = \alpha^{'}(\pi(\gamma))(x) - x
= \alpha(\gamma)(x) - x.$$
Since this contradicts the fact that $c$ is not a virtual coboundary,
we conclude that $c^{'}$ is also not a virtual coboundary. Hence
$H^{1}_{v}(X,\alpha^{'})\ne 0$.$\Box$

\medskip
Our next result shows that for Notherian actions of virtually
nilpotent groups on zero-dimensional compact abelian groups, the
virtual first cohomology group vanishes. In view of Theorem \ref{ns} and
Lemma \ref{simple} this completes the proof of Theorem \ref{nil}.
   
\medskip
\begin{theorem}\label{cohom}
Let $\Gamma$ be  a finitely
generated virtually nilpotent group, and let $(X,\alpha)$ be a 
Noetherian action of $\Gamma$ on  a  zero-dimensional group $X$.
Then  $H^{1}_{v}(\alpha)= \{ 0\}$
\end{theorem}
{\it Proof. \/}
Let $\Lambda\subset\Gamma$ be a finite index
subgroup, and let $c :\Lambda\rightarrow X$ be a 
virtual 1-cocycle of $\alpha$. Suppose
$\Gamma_0\subset \Gamma$ is a nilpotent subgroup of finite index.
Then $c$ is equivalent to $c|_{\Gamma_0\cap \Lambda}$. Hence, to show that $c$ is
a virtual coboundary, without loss of generality we may assume that 
$\Lambda\subset \Gamma_0$. 

We will use induction on $h(\Lambda)$, the Hirsch number of $\Lambda$.
If $h(\Lambda)) = 0$ then $\Lambda$ is finite, and hence $c$ is a virtual
coboundary. Suppose $h(\Lambda) \ge 1$. Then $\Lambda$ is a finitely
generated infinite nilpotent group. Hence the center of $\Lambda$ is
also finitely generated and infinite (\cite{Ho},
Proposition 2.8). We
choose an element $\gamma_0$ in the center of $\Lambda$ that has
infinite order. For $m\ge 1$, we define $K_m = (\alpha(\gamma_{0}^{m}) -I)(X)$.  
 Since $\alpha(\gamma_0)$ commutes with $\alpha(\gamma)$ for all
$\gamma\in\Gamma$, it follows that each $K_m$ is a closed $\alpha$-invariant
subgroup. As $\alpha$ is Noetherian, the collection $\{ K_m : m\ge
1\}$ has a minimal element $K$. Clearly $K$ is of the form  
$(\alpha(\gamma_1) - I)(X)$, where $\gamma_1 = \gamma_{0}^{n}$ for some $n\ge 1$.
Let $\beta$ denote the induced algebraic $\Gamma$-action on
the quotient $X/K$. If $P$ denotes the projection map from $X$ to
$X/K$ then $P\circ c$ is virtual 1-cocycle of the action $\beta$.
We note that the $\beta(\gamma) = I$ for all 
$\gamma$ in the infinite cyclic subgroup generated by $\gamma_1$. By
Lemma \ref{lem}
and the induction hypothesis we deduce that $P\circ c$ is a virtual
coboundary of $\beta$. We find a finite index subgroup
$\Lambda_1\subset \Lambda$ and $p\in X/K$ such that 
$$ P\circ c(\gamma) = \beta(\gamma)(p) - p \ \forall \gamma\in\Lambda_1.$$ 
 We choose $q\in X$ such that
$P(q) = p$. Let $c_1 :\Lambda_1\rightarrow X$ denote the virtual
1-cocycle defined by $c_1 = c - c_q$. Since $P$ is an equivariant
map it follows that $P\circ c_1 = 0$, i.e., the image of $c_1$ is
contained in $K = (\alpha(\gamma_1) - I)(X)$. We choose $l\ge 1$ such
that $\gamma_{1}^{l}\in\Lambda_1$. From the minimality of $K$ we
deduce that $K = (\alpha(\gamma_{2}) - I)(X)$, where $\gamma_2 =
\gamma_{1}^{l}$. We find 
$x\in X$ such
that $\alpha(\gamma_2)(x) - x = c_1(\gamma_2)$. Let $c_2 : \Lambda_1
\rightarrow X$ be the virtual 1-cocycle defined by 
$c_2(\gamma) = c_1(\gamma) - c_x(\gamma)$. Then $c_2$ is cohomologous
to $c_1$ and $c_2(\gamma_2) = 0$.
 Let $F\subset X$ denote the set of points that are fixed by
 $\alpha(\gamma_2)$. Since $\alpha(\gamma_2)$ commutes with
 $\alpha(\gamma)$ for all $\gamma$, $F$ is a 
$\alpha$-invariant closed subgroup. We note that for any $\gamma\in
\Lambda_1$, $c_2(\gamma_2\gamma) = c_2(\gamma_2) +
\alpha(\gamma_2)(c_2(\gamma)) = \alpha(\gamma_2)(c_2(\gamma))$.
As $\gamma_2$ lies in the center of $\Gamma$ we also
have,
$$ c_2(\gamma_2\gamma) = c_2(\gamma\gamma_2) = c_2(\gamma) +
\alpha(\gamma)(c_2(\gamma_2)) = c_2(\gamma).$$
This shows that for all $\gamma\in\Lambda_1$ the element $c_2(\gamma)$
lies in $F$.  
Let $\alpha_F$ denote the restriction of $\alpha$ to $F$. Then $c_2$
can be viewed as a virtual 1-cocycle of $\alpha_F$. Since
$\alpha(\gamma)|_{F} = I$ for all $\gamma$ in the infinite cyclic
group generated by $\gamma_2$, from the previous lemma
 and the induction hypothesis
we conclude that $c_2$ is a virtual coboundary. 
Since $c$ is cohomologous to
$c_2$ this completes the proof.
$\Box$

\medskip
We conclude this paper with the following question :
{\q \ }Does Theorem 1.1 hold if $\Gamma$ is
 polycyclic-by-finite ?

\medskip
In view of Theorem \ref{ns} and Lemma \ref{simple} this is true if
and only if Theorem \ref{cohom} holds for polycyclic-by-finite groups.

\end{document}